\documentclass[12pt,reqno,intlim]{amsart}

\usepackage{amssymb,amsmath}

\usepackage[arrow,matrix,curve]{xy}

\newdir{ >}{{}*!/-5pt/\dir{>}}

\newcommand{\nc}{\newcommand}

\nc{\bC}{\bold{C}} \nc{\bN}{\Bbb{N}} \nc{\cF}{\mathcal{F}}
\nc{\cE}{\mathcal{E}} \nc{\cR}{\mathcal{R}} \nc{\cM}{\mathcal{M}}
\nc{\al}{\alpha} \nc{\bt}{\beta} \nc{\gm}{\gamma} \nc{\dl}{\delta}
\nc{\om}{\omega} \nc{\sg}{\sigma} \nc{\Sg}{\Sigma} \nc{\vf}{\varphi}
\nc{\ve}{\varepsilon} \nc{\os}{\overset} \nc{\ol}{\overline}
\nc{\ul}{\underline} \nc{\us}{\underset} \nc{\sbs}{\subset}
\nc{\bsl}{\backslash} \nc{\Ra}{\Rightarrow}
\nc{\lra}{\longrightarrow} \nc{\all}{\allowdisplaybreaks}

\nc{\Codes}{\operatorname{{\bold{Codes}}}}
\nc{\RegMono}{\operatorname{\mathcal{R}{\rm{eg}\mathcal{M}{\rm{ono}\!}}}}
\nc{\RegEpi}{\operatorname{\mathcal{R}{\rm{eg}\mathcal{E}{\rm{pi}\!}}}}
\nc{\Mn}{\operatorname{\mathcal{M}{\rm{ono}\!}}}
\nc{\Ep}{\operatorname{\mathcal{E}{\rm{pi}\!}}}
\nc{\Rg}{\operatorname{\mathcal{R}{\rm{eg}\!}}}
\nc{\Ob}{\operatorname{Ob\!}}

\numberwithin{equation}{section}

\newtheorem{theo}{\ \ \ Theorem}[section]

\newtheorem{prop}[theo]{\ \ \ Proposition}

\theoremstyle{definition}
\newtheorem{exmp}[theo]{\ \ \ Example}

\theoremstyle{remark}

\begin{document}
\title[]
{Sufficient conditions for solvability of linear Diophantine equations, and Frobenius numbers}

\author{Eteri Samsonadze}

\maketitle

\begin{abstract}
The sufficient conditions for solvability  of a linear Diophantine equation 
$\sum_{i=1}^{n}a_ix_i=b$ (with $a_1,a_2,...,a_n\in \mathbb{N}$) in non-negative integers $x_1,x_2,...,x_n$ are given. The explicit formulas are given for Frobenius numbers $g(a_1,a_2,...,a_n)$, for some particular cases,. Besides, a new recurrent method of studying the problem of solvability of a linear Diophantine equation in non-negative integers is proposed. This recurrent method is used for the problem of finding Frobenius numbers $g(a_1,a_2,...,a_n)$ for any $n\geq 3$; the example is given for the case $n=5$.

\bigskip


\end{abstract}
\vskip+10mm

\section{Introduction}
A linear Diophantine equation
\begin{equation} \label{1.1.1}
\sum_{i=1}^{n}a_ix_i=b
\end{equation}
is considered, where $n\geq 2$, $a_1,a_2,...,a_n$ are coprime numbers; $b$ and $x_i$ $(i=1,2,...,n)$ are non-negative integers. 

As is well-known \cite{B},  if  $(a_1,a_2)=1$ and 
\begin{equation}
b>a_1a_2-(a_1+a_2), \label{1..1}
\end{equation}
 then the equation
$$a_1x_1+a_2x_2=b$$
\noindent is solvable in non-negative integers $x_1$ and $x_2$. This implies that, under condition (\ref{1..1}),  equation (\ref{1.1.1}) also is solvable in non-negative integers if $(a_1,a_2)=1$. In the present paper, we consider the general case where $a_1, a_2,...,a_n$ are coprime numbers, and prove that if 
$$\left[\frac{b}{M}\right]\geq \left[n-\frac{\sum_{i=1}^{n}a_i+r}{M}\right],$$ then equation (\ref{1.1.1}) is solvable in non-negative integers. Here $M$ is the least common multiple of the numbers $a_1,a_2,...,a_n$; $r$ is the remainder of $b$ modulo $M$; the symbol $\left[ x \right]$ denotes the largest integer not exceeding a number $x$. This fact, in particular,  implies that equation (\ref{1.1.1}) is solvable in non-negative integers if  $(a_1,a_2,...,a_n)=1$ and either 
$$b\geq nM-\sum_{i=1}^{n}a_i$$
or 
$$b\geq (n-1)M.$$ Moreover, if $(a_1,a_2,...,a_n)=1$ and
$$\sum_{i=1}^na_i\geq (n-2)M+2$$
 (or equivalently $$(n-1)M-\sum_{i=1}^{n}a_i\leq M-2),$$
 then equation (\ref{1.1.1}) is solvable in non-negative integers provided that
 $$b>(n-1)M-\sum_{i=1}^{n}a_i.$$  

In the paper, we also show that if $(a_1,a_2)=d \neq 1$, then equation (\ref{1.1.1}) is solvable in non-negative integers for $$b=(a'_1a'_2-(a'_1+a'_2)+k)d,$$ where $a'_i=\frac{a_i}{d}$ $(i=1,2)$ and $k$ is any natural number. 

The paper deals also with the problem of finding Frobenius numbers $g(a_1,a_2,...,a_n)$ for some particular cases. It is proved that $$g(a_1,a_2,...,a_n)=a_1-1$$ if $$3\leq a_1\leq a_2\leq ...a_n$$ and any number from the set $\lbrace a_1,a_1+1,a_1+2,...,a_1^2-a_1-1\rbrace$ is equal to at least one of $a_1,a_2,...,a_n$. Further, it is proved that if $k,a_2,a_3,...,a_n$ are natural numbers,  $(2k,a_2k,a_3k,...,a_{n-1}k,a_n)=1$ and $c$ is the smallest odd number among the numbers $a_2,a_3,...,a_n$, then $$g(2k,a_2k,a_3k,...,a_{n-1}k,a_n)=k(c-2)+a_n(k-1).$$ In particular, 
$$g(2k,ak,b)=k(c-2)+b(k-1),$$
where $k,a,b$ are natural numbers, $(2k,ak,b)=1$; $c=b$ if $a$ is an even number, and $c=min(a,b)$ otherwise.

Besides, a new recurrent method of studying the problem of solvability of a linear Diophantine equation (\ref{1.1.1}) in non-negative integers $x_1, x_2,..., x_n$ is proposed. In contrast to the well-known recurrent methods, this method reduces the case $b=b_0$ to the cases where $b= \left[\frac{b_0}{2}\right]-k_i$, where $k_i$ are some integer non-negative numbers. This method enables us to find Frobenius numbers $g(a_1,a_2,...,a_n)$ for any natural $n\geq 3$. In the paper, applying this method, the example of calculating a Frobenius number $g(a_1,a_2,...,a_n)$ is given for $n=5$.

\section{On Diaphantine linear equations and Frobenius numbers}

As is shown in \cite{S}, for the number $P(b)$ of non-negative integer solutions $(x_1,x_2,...,x_n)$ of the equation 
\begin{equation} \label{2.2}
\sum_{i=1}^{n}a_ix_i=b,
\end{equation} where $n\geq 2$, $a_1,a_2,...,a_n$ are coprime natural numbers and $b$ is a non-negative integer, the following equality is valid:
\begin{equation} \label{2.3}
P(b)=\sum_{k=0}^{s}l_k \overline{C}^{n-1}_{b'+n-1-k}
\end{equation}
\noindent where $$b'=\left[\frac{b}{M}\right];$$ $M$ is the least common multiple of the numbers $a_1,a_2,...,a_n$; $l_k$ $(k=0,1,...,s)$ is the number of integer solutions of the system
$$\sum_{i=1}^{n}a_ix_i=r+Mk, ~0\leq x_i\leq \frac{M}{a_i}-1 ~ (i=1,2,...,n);$$
$r$ denotes the remainder of $b$ modulo $M$; $\overline{C}_m^{t}$ denotes $C_m^{t}$ if $m\geq t$ and denotes $0$ otherwise (here $C_m^{t}$ is the binomial coefficient ${\begin{pmatrix}m\\t\\\end{pmatrix}}$);
\begin{equation} \label{2.4}
s=\left[n-\frac{\sum_{i=1}^{n}a_i+r}{M}\right],
\end{equation}

$$(0\leq s\leq n-1).$$
According to \cite{S}, for the case of coprime $a_1,a_2,...,a_n$, we have
\begin{equation} \label{2.6}
\sum_{i=1}^{s}l_i=\frac{M^{n-1}}{a_1a_2...a_n}.
\end{equation}

Since $l_i\geq 0$ ($i=0,1,...,s$), formula (\ref{2.6}) implies that at least one of $l_0,l_1,...,l_s$ is different from zero. Formula (\ref{2.3}) implies that equation (\ref{2.2}) is solvable if and only if there is a non-zero number $l_k$ with $k\leq b'$ ($0\leq k\leq s$). This, in particular, implies that if $l_0\neq 0$ (or equivalently $P(r)\neq 0$), then equation (\ref{2.2}) is solvable.

\begin{theo}
If $(a_1,a_2,...,a_n)=1$ and
\begin{equation}
\left[\frac{b}{M}\right]\geq \left[n-\frac{\sum_{i=1}^{n}a_i+r}{M}\right],
\end{equation}
then equation (\ref{2.2}) is solvable in non-negative integers.
\end{theo}

\begin{proof}
If $b'\geq s$, then $\overline{C}_{b'+n-1-k}^{n-1}\neq 0$ for $k=0,1,...,s$. Since at least one of the numbers $l_0,l_1,...,l_s$ is different from zero, formula (\ref{2.3}) implies that $P(b)\neq 0$.
\end{proof}

Theorem 2.1 implies 
\begin{prop}
If $(a_1,a_2,...,a_n)=1$ and either
$$b\geq nM-\sum_{i=1}^{n}a_i$$
 or 
$$b\geq (n-1)M,$$
then equation (\ref{2.2}) is solvable in non-negative integers.
\end{prop}




Besides, according to \cite{S}, we have
\begin{prop}
Let $(a_1,a_2,...,a_n)=1$. If 
$$\sum_{i=1}^{n}a_i\geq (n-2)M-2$$
(or, equivalently, $(n-1)M-\sum_{i=1}^{n}a_i\leq M-2$), then equation (2.1) is solvable in non-negative integers for
\begin{equation} \label{2.10.10}
b>(n-1)M-\sum_{i=1}^{n}a_i.
\end{equation}
\end{prop}

\begin{exmp}
The equation
\begin{equation} \label{2.9}
6x_1+70x_2+105x_3=b
\end{equation}
is solvable in non-negative integers for $b>239$. Indeed, here $\sum_{i=1}^{3}a_i=181$, $M=210$. In this case, we have 
$$s=\left[3-\frac{181+r}{210}\right],$$ 
where $r$ is the remainder of $b$ modulo $210$. 

Hence $s=1$ if $r>29$; and $s=2$ if $r\leq 29$. Therefore, Theorem 2.1 implies that equation (\ref{2.9}) is solvable if $b>210$ and $r>29$; it is also solvable if $b\geq 2\cdot 210$ and $r\leq 29$. Thus, equation (2.7) is solvable if $b>210+29=239$.
\end{exmp}


\vskip+2mm
Theorem 2.1 implies the well-known fact \cite{B} that if $(a_1,a_2)=1$, then the equation 
\begin{equation}  \label{2.13}
a_1x_1+a_2x_2=b
\end{equation}
is solvable in non-negative integers if
\begin{equation}  \label{2.14.14}
b>a_1a_2-(a_1+a_2).
\end{equation}
 Indeed, here $$s=\left[2-\frac{a_1+a_2+r}{a_1a_2}\right],$$ where $r$ is the remainder of $b$ modulo $a_1a_2$. Hence $s=0$ if $r>a_1a_2-(a_1+a_2)$; and $s=1$ if $r\leq a_1a_2-(a_1+a_2)$. Therefore, Theorem 2.1 implies that equation (\ref{2.13}) is solvable if $b> 0$ and $r> a_1a_2-(a_1+a_2)$; it is solvable also in the case where $b\geq a_1a_2$ and $r\in \left[0,a_1a_2-(a_1+a_2)\right]$. Thus, equation (\ref{2.13}) is solvable for any $b$ satisfying condition (\ref{2.14.14}). This implies

\begin{prop}
The equation
$$\sum_{i=1}^{n}a_ix_i=b$$
 is solvable in integer non-negative numbers $x_1,x_2,...,x_n$ if, among the numbers $a_1,a_2,...,a_n$, there exist coprime numbers $a_i$ and $a_j$ with 
 $$b>a_ia_j-(a_i+a_j).$$
\end{prop}

\begin{proof}
Assume that $(a_1,a_2)=1$. Then, as it was noted above, equation (\ref{2.13}) is solvable if the condition (\ref{2.14.14}) is satisfied. If $(y_1,y_2)$ is the solution of this equation, then the $n$-tuple $(y_1,y_2,0,0,...,0)$ is obviously a solution of the original equation.
\end{proof}

Recall that an integer $b_0$ is called the Frobenius number (and is usually denoted by the symbol $g(a_1,a_2,...,a_n)$) if $b=b_0$ is the greatest number such that the equation  
$$\sum_{i=1}^{n}a_ix_i=b$$
has no solution in integer non-negative numbers $x_1,x_2,...,x_n$.  

The problem of finding Frobenius numbers is closely related to some problems from various fields of mathematics. A lot of papers are devoted to this issue. However, the explicit  formulas for calculating Frobenius numbers are found only for the cases $n=2$  and $n=3$. Namely, it is proved that 
\begin{equation}
g(a_1,a_2)=a_1a_2-(a_1+a_2)
\end{equation}
if $(a_1,a_2)=1$ (see, e. g., \cite{B}, \cite{BS}, \cite{Sy}). For the case where $n=3$,  various formulas for calculating Frobenius numbers $g(a_1,a_2,a_3)$ are given in \cite{T2}; they correspond to various conditions imposed on $a_1,a_2,a_3$.
\vskip+3mm
Proposition 2.5 implies 
\begin{prop}
Let $a_1,a_2,...,a_n$ be natural numbers. If, among them, there exist coprime numbers $a_i$ and $a_j$, then
$$g(a_1,a_2,...,a_n)\leq a_ia_j-(a_i+a_j).$$
\end{prop}

\begin{prop}
Let $a_1,a_2,...,a_n$ be natural numbers. If, among them, there exist coprime numbers $a_i$ and $a_j$, then
$$g(a_1,a_2,...,a_n)<M,$$
where $M$ is the least common multiple of $a_1,a_2,...,a_n$.
\end{prop}

As is well-known \cite{T2}, 
$$g(a_1,a_2,a_3)=g(a_1,a_2)$$
if $a_3>g(a_1,a_2)$. This result can be generalized. Namely, we have the following statement.
\begin{prop} Let $k$ be a natural number with $1\leq k< n$. If $a_j>g(a_1,a_2,...,a_k)$ for any $j$ with $k<j\leq n$, then
$$g(a_1,a_2,...,a_k,a_{k+1},...,a_n)=g(a_1,a_2,...,a_k).$$
\end{prop}

\begin{proof}
The equation
$$a_1x_1+a_2x_2+...+a_kx_k+a_{k+1}x_{k+1}+...+a_nx_n=g(a_1,a_2,...,a_k)$$
has no solution in non-negative integers. Indeed, otherwise one would have $x_j=0$ for any $j>k$, in a solution $(x_1,x_2,...,x_n)$, as it follows from the assumption of the proposition; this would imply that the equation
$$a_1x_1+a_2x_2+...+a_kx_k=g(a_1,a_2,...,a_k).$$
is solvable that contradicts the definition of a Frobenius number. 

If $b>g(a_1,a_2,...,a_k)$, the equation 
$$\sum_{i=1}^{k}a_ix_i=b$$
is solvable as it follows from the definition of a Frobenius number. Let $(y_1,y_2,...,y_k)$ be its solution. Then obviously the $n$-tuple $$(y_1,y_2,...,y_k,0,0,...,0)$$ is a solution of the equation
$$\sum_{i=1}^{n}a_ix_i=b.$$
Thus, this equation is solvable for $b>g(a_1,a_2,...,a_k)$ and is not solvable for $b=g(a_1,a_2,...,a_k)$. Therefore, $$g(a_1,a_2,...,a_k, a_{k+1},...,a_n)=g(a_1,a_2,...,a_k).$$
\end{proof}

\begin{exmp}
Prove that 
\begin{equation} \label{2.16}
g(6,70,105,317)=239.
\end{equation}
The equation
$$6x_1+70x_2+105x_3=b$$
is solvable in non-negative integers for $b>239$ (see Example 2.4), while the equation
$$6x_1+70x_2+105x_3=239$$
is not solvable. Indeed, assume that the latter equation has a solution $(x_1,x_2,x_3)$. It is easy to observe that $x_3$ in it should be odd and less than 3. But if $x_3=1$, then this implies that $$6x_1+70x_2=134,$$ or, equivalently, 
$$3x_1+35x_2=67.$$
However, this equation is not solvable since
$$67=3\cdot 35-(3+35)=g(3,35).$$
Thus, $$g(6,70,105)=239;$$
and since $317>g(6,70,105)$, Proposition 2.8 implies equality (\ref{2.16}).
\end{exmp}
\vskip+5mm

\begin{prop}
If $(2,a_2,a_3,...,a_n)=1$, then
\begin{equation}
g(2,a_2,a_3,...,a_n)=c-2,
\end{equation}
where $c$ is the smallest odd number among the numbers $a_2,...,a_n$.
\end{prop}

\begin{proof}
First note that there is an odd number, among the numbers $a_2,...,a_n$, since $(2,a_2,a_3,...,a_n)=1$. 

Proposition 2.5 implies that the equation
\begin{equation}
2x_1+a_2x_2+a_3x_3+...+a_nx_n=b
\end{equation}
is solvable in integer non-negative numbers for $b>c-2$ (since $2c-(2+c)=c-2$).

Consider the equation
$$2x_1+a_2x_2+a_3x_3+...+a_nx_n=c-2.$$
If one removes the summands, where the variable's coefficients are greater than $(c-2)$, from the left part of this equation, then obviously the number of integer non-negative solutions of this equation will not be changed. But the obtained equation will be insolvable. Indeed, no summands with odd coefficients will remain on the left-hand part of the equation, while one will have an odd number on its right-hand part.

Thus, equation (2.13) is insolvable for $b=c-2$, and is solvable for $b>c-2$. This implies (2.12).
\end{proof}

\begin{prop}
Let $k$ and $a_2,...,a_n$ be natural numbers. If $$(2k,a_2k,a_3k,...,a_{n-1}k,a_n)=1$$ and $c$ is the smallest odd number among $a_2,a_3,...,a_n$, then
\begin{equation}
g(2k,a_2k,...,a_{n-1}k,a_n)=k(c-2)+a_n(k-1).
\end{equation}
In particular, 
\begin{equation}
g(2k,ak,b)=k(c-2)+b(k-1),
\end{equation}
where $k,a,b$ are natural numbers, $(2k,ak,b)=1$; $c=b$ if $a$ is even, and $c=min(a,b)$ otherwise.
\end{prop}

\begin{proof}
If $(a_1,a_2,a_3,...,a_n)=1$, $a_i=ka_i'$ $(i=1,2,...,n-1)$, then 
$$g(a_1,a_2,...,a_{n-1},a_n)=kg(a_1',a'_2,...,a'_{n-1},a_n)+a_n(k-1),$$
as it follows from \cite{J}, \cite{BS}. This implies that
$$g(2k,a_2k,...,a_{n-1}k,a_n)=kg(2,a_2,a_3...,a_{n-1},a_n)+a_n(k-1).$$
Here $(2,a_2,...,a_{n-1},a_n)=1$ since $(2k,a_2k,...,a_{n-1}k,a_n)=1.$ Therefore, formula (2.12) implies (2.14).
\end{proof}




\begin{exmp}
Formulas (2.14) and (2.15) imply that 
\vskip+1mm
(a) $g(40,80,100,101)=1979$
(since here $k=20$, $c=\frac{100}{20}=5$, $a_n=101$, and $20\cdot (5-2)+101\cdot 19=1979$).
\vskip+2mm
(b) $g(22,110,115)=2393$ (since $g(22,110,115)=g(2\cdot 11, 10\cdot 11,  115)=11\cdot (115-2)+115\cdot 10=2393$).
\end{exmp}


\begin{prop}
We have
\begin{equation}
g(a,a+1,a+2,...,g(a,a+1))=a-1
\end{equation}
for $a\geq 3$.
\end{prop}

\begin{proof}
First note that since 
\begin{equation}
g(a,a+1)=a^2-a-1,
\end{equation}
we have
\begin{equation}
g(a,a+1)>a+1
\end{equation}
for $a\geq 3$.

It is obvious that the equation 
\begin{equation}
ax_1+(a+1)x_2+(a+2)x_3+...+g(a,a+1)x_n=b
\end{equation}
is solvable in integer non-negative numbers for $b=a, a+1, a+2,...,g(a,a+1)$. Moreover, since $(a,a+1)=1$, Proposition 2.5 implies that it is solvable for $b>g(a,a+1)$. Thus, this equation is solvable for $b\geq a$. It is obvious that this equation is not solvable for $b=a-1$. This implies (2.16).
\end{proof}

Proposition 2.13 implies the following statement.

\begin{prop} For any natural numbers $k$, $a\geq 3$ and $c_1,c_2,...,c_k\geq a$, we have
\begin{equation}
g(a,a+1,a+2,...,g(a,a+1),c_1,c_2,...,c_k)=a-1
\end{equation}

\end{prop}

\begin{proof} If equation (2.19) is solvable in integer non-negative numbers $x_1,x_2,...,x_n$ for $b\geq a$, then the equation
$$ax_1+(a+1)x_2+(a+2)x_3+...+g(a,a+1)x_n+c_1y_1+c_2y_2+...+c_ky_k=b$$
also is solvable in integer non-negative numbers $x_1,x_2,...,x_n, y_1,y_2,...,y_k$ for $b\geq a$. For $b=a-1$, the latter equation is obviously insolvable. This implies (2.20).
\end{proof}

Formulas (2.17), (2.18), (2.20) imply

\begin{theo}
If $3\leq a_1\leq a_2\leq ...\leq a_n$, and all numbers from the set $\lbrace a_1,a_1+1,a_1+2,...,a_1^2-a_1-1\rbrace$ present in the sequence $a_1,a_2,...,a_n$, then
$$g(a_1,a_2,...,a_n)=a_1-1.$$
\end{theo}

\begin{exmp}
We have
$$g(4,5,6,6,7,8,9,10,11,12)=3;$$
(here $a_1=4$, $a_1^2-a_1-1=11$).
\end{exmp}

Theorem 2.15 implies

\begin{prop}
If $a\geq 3$ and $n\geq a^2-2a-1$, then
\begin{equation}
g(a,a+1,a+2,...,a+n)=a-1
\end{equation}

\end{prop}

\begin{proof}
It suffices to observe that if $n\geq a^2-2a-1$, then $a+n\geq a^2-a-1$.
\end{proof}

Thus, for sufficiently large $n$, we have formula (2.21).
\vskip+2mm
Note that formula (2.21) follows also from the following formula from \cite{Se}:
$$g(a,ha+d,ha+2d,...,ha+nd)=ha\left[\frac{a-2}{n}\right]+(h-1)a+d(a-1),$$
where $d$ and $h$ are natural numbers and $(a,d)=1$.
Indeed, if $a\geq 3$, then $(a-2)<a^2-2a-1$. Hence $\left[\frac{a-2}{n}\right]=0$ if $n\geq a^2-2a-1$. This formula, for $h=d=1$ implies (2.21).

Note that formulas for Frobenius number $g(a_1,a_2,...,a_n)$, for the case where the sequence $a_1,a_2,...,a_n$ is an arithmetic progression, are given, for instance, in \cite{B}, \cite{R}, \cite{T}.

\begin{theo}
If $(a_1,a_2)=d\neq 1$, then the equation
\begin{equation}  \label{2.18}
\sum_{i=1}^{n}a_ix_i=b
\end{equation}
is solvable in non-negative integers for 
\begin{equation}  \label{2.19}
b=(a'_1a'_2-(a'_1+a'_2)+k)d,
\end{equation} where $$a'_i=\frac{a_i}{d}$$ $(i=1,2)$, and $k$ is an arbitrary natural number. In particular, if $d=2$, then this equation is solvable for any even integer $b$ with $$b>2(a'_1a'_2-(a'_1+a'_2)).$$
\end{theo}

\begin{proof}
Consider the equation
$$a'_1x_1+a'_2x_2=b.$$
We have $(a'_1,a'_2)=1$. Therefore, this equation is solvable for $b$ satisfying the inequality
$$b>a'_1a'_2-(a'_1+a'_2),$$ i.e., for $b$ that is given by
$$b=a'_1a'_2-(a'_1+a'_2)+k,$$ where $k$ is any natural number. Therefore, the equation
$$a_1x_1+a_2x_2=b$$
is solvable for any $b$ that is given by (\ref{2.19}). This implies that if $n\geq 2$, then also equation (\ref{2.18}) is solvable for this $b$.
\end{proof}

Recall the recurrent formula for the number $P(b)$ of non-negative integer solutions of a linear Diophantine equation obtained in \cite{S}. 

\begin{theo} \cite{S}
Let $P(b)$ be the number of non-negative integer solutions of the equation 
$$\sum_{i=1}^{n}a_ix_i=b$$
where $n\geq 2$, $a_1,a_2,...,a_n$ are coprime numbers, $b$ is a non-negative integer. Then
$$P(b)=\sum_{k=0}^{s_1}P^{\ast}(2k)P(b'-k)$$
for an even $b$, and
$$P(b)=\sum_{k=0}^{s_2}P^{\ast}(2k+1)P(b'-k)$$
for an odd $b$, where $$b'=\left[\frac{b}{2}\right], \;\;\; s_1=\left[\frac{\sum_{i=1}^{n}a_i}{2}\right], \;\;\;  s_2=\left[\frac{\sum_{i=1}^{n}a_i-1}{2}\right],$$ while
$P^{\ast}(d)$ is the number of integer solutions of the system
$$\sum_{i=1}^{n}a_ix_i=d, \;\;\;\; 0\leq x_i\leq 1,\;\;\; (i=1,2,...,n).$$
\end{theo}

Theorem 2.19 implies the following necessary and sufficient condition for the solvability of a linear Diophantine equation in non-negative integers.

\begin{prop}
 Let $n\geq 2$, $a_1,a_2,...,a_n$ be coprime numbers.
 
 (i) If $b$ is an even integer, then the equation
\begin{equation} \label{2.20}
\sum_{i=1}^{n}a_ix_i=b
\end{equation} 
 is solvable in non-negative integers if and only if, among the numbers $0,1,2,...,s_1$, there is a number $k$ such that $$P^{\ast}(2k)\neq 0 \;\;\; and \;\;\; P(b'-k)\neq 0;$$
here $b'=\frac{b}{2}$.

(ii) If $b$ is an odd integer, then equation (\ref{2.20}) is solvable in non-negative integers if and only if, among the numbers $0,1,2,...,s_2$, there is a number $k$ such that $$P^{\ast}(2k+1)\neq 0 \;\;\; and \;\;\; P(b'-k)\neq 0;$$ 
here $b'=\left[\frac{b}{2}\right]$.
\end{prop}
\vskip+3mm
Let $a_1,a_2,...,a_n\in\mathbb{N}$, and $T$ be the set of all numbers that can be represented in the form
$$\sum_{i=1}^{n}a_ix_i$$ 
with $x_i$ being equal to either $0$ or $1$ ($i=1, 2, ..., n)$, i.e.:
$$T=\lbrace \sum_{i=1}^{n}a_ix_i: x_i=0 \; or \; 1  (i=1, 2, ..., n)\rbrace.$$

Let $D=\lbrace d_1,d_2,...,d_l\rbrace$ be the subset of $T$ containing only even integers, while $C=\lbrace c_1,c_2,...,c_m\rbrace$ be the subset of $T$ containing only odd integers. Let $D'=\lbrace d'_1,d'_2,...,d'_l\rbrace$, where $d'_i=\frac{d_i}{2}$ $(i=1,2,...,l)$, while $C'=\lbrace c'_1,c'_2,...,c'_m\rbrace$, where $c'_j=\left[\frac{c_j}{2}\right]$ $(j=1,2,...,m)$. 
\vskip+3mm
It is obvious that if $b$ is a non-negative integer and $b$ is divisible by some element of the set $T$, then equation (2.24) is solvable. 
\vskip+1mm
Proposition 2.20 implies the following statement.
\begin{theo}
Consider the equation
\begin{equation}  \label{2.21}
\sum_{i=1}^{n}a_ix_i=b
\end{equation} 
where $n\geq 2$, $a_1,a_2,...,a_n$ are coprime numbers, and $b$ is a non-negative integer. Let $b'=\left[\frac{b}{2}\right]$.

(i) Let $b$ be an even integer. The equation (\ref{2.21}) is solvable in non-negative integers if and only if at least one of the numbers $P(b'-d'_i)$ $(i=1,2,...,l)$ is non-zero.

(ii) Let $b$ be an odd integer. The equation (\ref{2.21}) is solvable in non-negative integers if and only if at least one of the numbers $P(b'-c'_j)$ $(j=1,2,...,m)$ is non-zero.   
\end{theo}

Theorem 2.21 enables us to study the solvability of a linear Diophantine equation in non-negative integers employing a recurrent method. 
\vskip+2mm
Let $a$ be the smallest number among the naturals $a_1, a_2,...,a_n$. It is obvious that if $b<2a$, then the inequality
$$0\leq \sum_{i=1}^{n}a_ix_i<2a$$
implies that $0\leq x_i\leq1$ ($i=1,2,...,n$). This implies
\begin{prop}
Let $a$ be the smallest number among naturals $a_1, a_2,...,a_n$.
 If $b<2a$ and $b$ does not belong to the set $T$, then equation (2.25) is not solvable in non-negative integers. 
\end{prop}


The above-mentioned recurrent method for studying solvability of a linear Diophantine equation can be employed for calculating Frobenius numbers $g(a_1,a_2,...,a_n)$ for any $n\geq 3$.
\vskip+2mm

\begin{exmp}
Find the Frobenius number
$$g(6,8,11,13,15).$$
Since $(6,11)=1$ and $6\cdot 11-(6+11)=49$, we have $P(b)\neq 0$ for $b>49$, where $P(b)$ is the number of non-negative integer solutions of the equation
\begin{equation}  \label{2.23}
6x_1+8x_2+11x_3+13x_4+15x_5=b.
\end{equation}
It is not hard to verify that
$$T=\lbrace 0, 6, 8, 11,13,14, 15,17,19,21,23,24,25,$$
$$26,27,28,29,30,32,34,36,38,39,40,42,45,47,53\rbrace,$$
$$D'=\lbrace 0,3,4,7,12,13,14,15,16,17,18,19,20,21\rbrace,$$
$$C'=\lbrace 5,6,7,8,9,10,11,12,13,14,19,22,23,26\rbrace.$$
Therefore, Theorem 2.21 implies that if $b$ is an even number and $b\leq 49$, then $P(b)\neq 0$ if and only if at least one of the numbers $P(b'), P(b'-3), P(b'-4), P(b'-7), P(b'-12), P(b'-13), P(b'-14), P(b'-15), P(b'-16), P(b'-17), P(b'-18), P(b'-19), P(b'-20), P(b'-21)$ is non-zero. Moreover, if $b$ is an odd number and $b\leq 49$, then $P(b)\neq 0$ if and only if at least one of the numbers $P(b'-5), P(b'-6), P(b'-7), P(b'-8), P(b'-9), P(b'-10), P(b'-11), P(b'-12), P(b'-13), P(b'-14), P(b'-19), P(b'-22), P(b'-23), P(b'-26)$ is non-zero. 

Moreover, $P(b)\neq 0$ if $b$ is divisible by at least one number from the set $T$; $P(b)=0$ if $b\in \lbrace 1,2,3,4,5,7,9,10\rbrace$, as it follows from Proposition 2.22. 

For an odd $b$, we obtain
$P(49)\neq 0$ since $\left[\frac{49}{2}\right]-5=19\in T$; $P(47)\neq 0$ since $47\in T$; $P(45)\neq 0$ since $45\in T$; $P(43)\neq 0$ since $\left[\frac{43}{2}\right]-8=13\in T$; $P(41)\neq 0$ since $\left[\frac{41}{2}\right]-6=14\in T$; $P(39)\neq 0$ since $39\in T$; $P(37)\neq 0$ since $\left[\frac{37}{2}\right]-5=13\in T$; $P(35)\neq 0$ since $\left[\frac{35}{2}\right]-9=8\in T$; $P(33)\neq 0$ since $\left[\frac{33}{2}\right]-5=11\in T$; $P(31)\neq 0$ since $\left[\frac{31}{2}\right]-9=6\in T$.

Since $29, 27, 25, 23, 21, 19, 17, 15, 13, 11\in T$, the numbers $P(29)$, $P(27)$, $P(25)$, $P(23), $P(21)$, $P(19)$,  $P(17), $P(15)$, $P(13)$, $P(11)$ are different from zero. 

As it was observed above, Proposition 2.22 implies that $P(10)=0$. 


Since $(6,8)=2$ and $3\cdot 4-(3+4)=5$, Theorem 2.18 implies that $P(b)\neq 0$ for an even $b$ with $b>2\cdot 5$, i.e., for $12, 14, 16...$.


Thus, for $b=10$, equation (\ref{2.23}) is not solvable, while for $b>10$ is solvable  (regardless of whether $b$ is even or odd). Therefore,
$$g(6,8,11,13,15)=10.$$
\end{exmp}
\vskip+5mm





\vskip+2mm
 
\textit{Author's addresses:}
 \textit{Eteri Samsonadze,
Retd., I. Javakhishvili Tbilisi State University,}
\textit{1 Tchavchavadze Av., Tbilisi, 0179, Georgia}, 

\textit{e-mail: etsamsonadze@gmail.com} 

\end{document}